\documentclass[10pt]{article}

\usepackage{setspace}
\usepackage[T2A]{fontenc}			
\usepackage[utf8]{inputenc}			
\usepackage[english]{babel}	

\usepackage{titleref}
\usepackage{hyperref}
\usepackage[lowercase]{theoremref}

\usepackage{graphicx}
\usepackage[font=small,labelsep=period]{caption}
\usepackage{subfigure}
\usepackage{float}
\usepackage{wrapfig}

\usepackage{titlesec}

\titleformat{\section}[block]{\Large\bfseries\filcenter}{\thesection.}{0.2em}{}

\usepackage{lipsum}
\usepackage{vmargin}
\setmarginsrb{2.5cm}{2.5cm}{2.5cm}{2.5cm}{0pt}{0mm}{0pt}{13mm}
\usepackage{icomma} 
\usepackage{enumitem}
\usepackage{euscript}
\usepackage{mathrsfs}

\usepackage{amsmath,amsfonts,amssymb,amsthm}					
\usepackage{mathtools}
\theoremstyle{plain}
\newtoks\thehProclaim
\newtheorem*{Proclaim}{\the\thehProclaim}
\newenvironment{proclaim}[1]
{\thehProclaim{#1}
\begin{Proclaim}}
{\end{Proclaim}}    
\newtheorem{lemma}{Lemma}
\newtheorem*{proposition}{Proposition}
\newtheorem*{remark}{Substitution remark}

\begin{document}

\title{On the~number of~tuples of~group elements satisfying a~first-order formula}
\author{Elena K. Brusyanskaya\footnote{The~work was supported by the~Russian Science Foundation, project no. 22-11-00075.}\\
\small \textsl{Faculty of~Mechanics and Mathematics of~Moscow State University}\\
\small \textsl{Moscow 119991, Leninskie gory, MSU.}\\
\small \textsl{Moscow Center of~Fundamental and Applied Mathematics}\\
 \small \textsl{ebrusianskaia@gmail.com}}

\maketitle

\begin{abstract}
Our result contains as special cases the~Frobenius theorem (1895) on the~number of~solutions to~the~equation $x^n=1$ in~a~finite group
and the~Solomon theorem (1969) on the~number of~solutions
in~a~group to systems of~equations with fewer equations than unknowns.
Instead of~systems of~equations, we consider arbitrary first-order formulae in~the~group language with constants.
Our~result substantially generalizes the~Klyachko--Mkrtchyan theorem (2014) on this topic.
\end{abstract}

\section{Introduction}
The~research was inspired by two classical results on divisibility in groups: the~Frobenius theorem (1895), and the~Solomon theorem (1969).

\begin{proclaim}{Frobenius theorem \cite{Frob95}}The~number of~solutions to the~equation $x^n=1$ in a~finite group $G$ is divisible by the~greatest common divisor of~the~group order and $n$ for any positive integer $n$.
\end{proclaim}
\begin{proclaim}{Solomon theorem \cite{Solo69}}  In any group, the~number of~solutions to a~system of~coefficient-free equations is divisible by the~order of~the~group provided the~number of~equations is less than the~number of~unknowns.
\end{proclaim}
These theorems have been generalized many times in various directions (see \cite{ACNT13}, \cite{AsTa01}, \cite{BrTh88},  \cite{BK22}, \cite{Frob03}, \cite{Hall36}, \cite{Isaa70}, \cite{Kula38}, \cite{Sehg62},  \cite{Yosh93}). For example, the~following generalization of~the~Solomon theorem was proved in \cite{GRV12}.

\begin{proclaim}{Gordon–Rodriguez-Villegas theorem~\cite{GRV12}} In any group, the~number of~solutions to a~finite system of~coefficient-free equations is divisible by the~order of~this group if the~rank of~the~matrix composed of~the~exponent sums of~the~$i$th unknown in the~$j$th equation is less than the~number of~unknowns. 
\end{proclaim}
There are non-obvious connections between the~Frobenius theorem and Solomon theorem. The~paper~\cite{BKV19} shows that these two theorems (and their generalizations, including the~Gordon–Rodriguez-Villegas theorem) are special cases of~a~general theorem called the \emph{divisibility theorem} (see below).

In this paper, we generalize classical divisibility theorems on systems of~equations to arbitrary first-order formulae. The~only known result in this direction is the~following theorem.
\begin{proclaim}{Klyachko--Mkrtchyan theorem (a simplified form) \cite{KM14}} 
	If the~rank of~the~matrix $A(\gamma)$ of~a~formula $\gamma$ is less than the~number of~free variables of~this formula, then the~number of~tuples of~group elements satisfying the~formula $\gamma$ is divisible by the~order of~centralizer of~the~set of~all coefficients of~$\gamma$. 
\end{proclaim}

  The~matrix of~a~formula is defined in the~next section. For now, let us confine ourselves to the~case of~the~quantifier-free formula $A(\gamma)=(a_{ij})$, where $a_{ij}$ is equal to the~exponent sum of~the~$j$th variable in the~$i$th atomic subformula.

We generalize the~Klyachko--Mkrtchyan theorem by eliminating~the~condition on the~rank of~the~formula matrix.

\begin{proclaim}{Main theorem}
	In any group, the~number of~tuples of~group elements satisfying a~formula $\gamma$ of~$m$ free variables is divisible by the~greatest common divisor of~the~centralizer of~the~set of~all coefficients of~$\gamma$ and the~integer  $n=\frac{\Delta_m}{\Delta_{m-1}}$, where $\Delta_i$ is the~greatest common divisor of~all minors of~order $i$ of~the~formula matrix. The~following conventions are assumed: $\Delta_i=0$ if $i$ is larger than the~number of~rows of~the~matrix; $\Delta_0=1$; $\frac{0}{0}=0$.
\end{proclaim}

We define the~greatest common divisor GCD$(G, n)$ of~a~group $G$ and an~integer $n$ as the~least common multiple of~the~orders of~all subgroups of~$G$ dividing $n$. The~divisibility is always understood in the~sense of~cardinal arithmetic: each infinite cardinal is divisible by all smaller nonzero cardinals (and zero is divisible by all cardinals and divides only zero). In the~case of~a~finite group $G$, we have GCD$(G, n)=$GCD$(|G|, n)$ by the~Sylow theorem and from the~fact that a~finite $p$-group contains subgroups of~all possible orders.
For example, the~matrix of~the~formula $\gamma(x_1, x_2)$ 

\begin{equation*}
\neg(x_1^3x_2=a)\vee(x_1x_2^3=b),
\end{equation*}
where $a, b \in G$ are coefficients, i.e. some fixed group elements,
has the~form

\begin{equation*}
A(\gamma) = \left(
\begin{array}{cc}
3 & 1\\
1 & 3
\end{array}
\right).
\end{equation*}

In this~example, there is one minor of~order 2 (equal to 8) and there are four minors of~order 1. We have GCD$(3, 1, 3, 1)=1$. Thus, the~theorem asserts that the~number of~tuples of~group elements satisfying $\gamma(x_1, x_2)$ is divisible by GCD$(C(a)\cap C(b), 8)$. Note that the~Klyachko–Mkrtchyan theorem has nothing to say about this case since the~rank of~the~matrix $A(\gamma)$ is equal to the~number of~free variables of~the~formula.

\emph{ Notation and conventions}. 
The~letter $\mathbb{Z}$ denotes the~set of~integers, GCD is the~greatest common divisor, and LCM is the~least common multiple. If $k \in \mathbb{Z}$, and $x$ and $y$ are elements of~a~group, then $x^y$, $x^{ky}$ and $x^{-y}$ denote $y^{-1}xy$, $y^{-1}x^ky$ and $y^{-1}x^{-1}y$, respectively. The~commutator $[x,y]$ is $x^{-1}y^{-1}xy$. If $X$ is a~subset of~a~group, then $|X|$, $\langle X\rangle$ and $C(X)$ are the~cardinality of~$X$, the~subgroup generated by $X$, and the~centralizer of~$X$. The~symbol $\langle g\rangle_n$ denotes the~cyclic group of~order $n$ generated by an~element $g$. The~free group of~rank $n$ is denoted by $F(x_1,\dots,x_n)$ or $F_n$. The~symbol $G*A$ denotes the~free product of~groups $G$ and $A$.

Let us note again that the~finiteness of~groups is not assumed by default; divisibility is always understood in the~sense of~cardinal arithmetic (an infinite cardinal is divisible by all nonzero cardinals not exceeding it), and GCD\,$(G, n) \overset{\text{def}}=$ LCM\,$({|H| \,\big| \,H \text{ is a~subgroup of~$G$, and $|H|$ divides $n$}})$.

The~word ``formula'' means a~first-order formula in the~group language with constants.

The~author thanks A. A. Klyachko and an~anonymous referee for valuable remarks.

\section{Main theorem}

Consider an~arbitrary first-order formula $\gamma$ in the~language  over a~group $G$. This language has two functional symbols $\cdot$ and $^{-1}$, and also, for each element of~$G$, there is a~constant symbol $g$.  

The~\emph{matrix of~a~formula $\gamma$} is an~integer matrix $A(\gamma)=(a_{ij})$, which is composed as follows. Each atomic subformula can be written in the~form
\begin{equation}\nonumber
    u=1,
\end{equation}
where $u\in G*F,$ and $F$ is the~free group generated by all (free and bound) variables of~$\gamma$. Thus, the~words~$u$ (possibly, different for different subformulae) can contain free and bound variables and elements of~$G$ (called \emph{coefficients} of~$\gamma$).

Now, we define the~labelled digraph $\Gamma(\gamma)$ of~the~formula $\gamma$ allowing loops and multiple edges. The~vertices of~ $\Gamma(\gamma)$ are all bound variables of~$\gamma$. Each atomic subformula containing bound variables has the~form
\begin{equation}\nonumber
    v_1(y_1)w_1(x_1,\dots, x_m)\dots v_r(y_r)w_r(x_1,\dots, x_m)=1,
\end{equation}

\noindent where $y_i$ are bound variables of~$\gamma$ (not necessarily different), $x_1,\dots, x_m$ are all (different) free variables of~$\gamma$, the~words $v_i(y_i)$ belong to the~free product $G*\langle y_i\rangle_\infty$ of~$G$ and the~infinite cyclic group generated by the~letter $y_i$, the~words $w_i(x_1,\dots, x_m)$ belong to the~free product $G*F(x_1,\dots, x_m)$ of~$G$ and the~free group with the~basis $x_1,\dots, x_m$. Let us connect $y_i$ and $y_{i+1}$ (subscripts modulo $r$) by a~directed edge labelled by an~integer tuple $(\alpha_1,\dots, \alpha_m)$, where $\alpha_j$ is the~exponent sum of~$x_j$ in $w_i$; loops labelled by zero tuples are excluded. We apply this construction to each atomic subformula containing bound variables.

For example, consider a~formula $\gamma(x_1,\; x_2)$:
\begin{equation*}
\forall y\; \exists z\; (ay(x_1x_2)^2yx_1b^zx_1^3=1)\wedge(b^z(x_1x_2)^2b^z=1)\vee\neg(x_1^2x_2=1) \text{, where $a, b$ are elements of~$G$.}
\end{equation*}

The~graph $\Gamma(\gamma)$ looks as follows:

\begin{figure}[!h]
    \centering
    \includegraphics[width=0.6\textwidth]{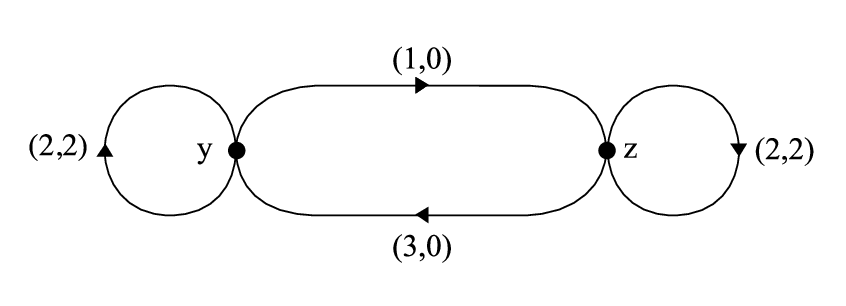}
    \caption{The~graph of~the~formula $\gamma(x_1,\;x_2)$}
    \label{fig:ex0}
\end{figure}

Now, in the~graph $\Gamma(\gamma)$, we choose cycles $c_1, c_2,\dots$ generating the~first homology group (e.g., generators of~the~fundamental groups of~all components can be taken), then we compose the~matrix $A(\gamma)$ of~$\gamma$ as follows: for each generating cycle $c_i$ we write a~row which is the~sum of~the~labels of~the~edges of~this cycle (the~edge labels are summed with signs plus or minus depending on the~orientation), then add rows, consisting of~exponent sums of~atomic subformulae containing no bound variables.

The~matrix $A(\gamma)$ depends on the~choice of~generating cycles, but the~integer linear hull of~its rows is uniquely determined by the~formula $\gamma$. The~matrix of~the~formula $\gamma(x_1, x_2)$ is

\begin{equation}\nonumber
A(\gamma) = \left(
\begin{array}{cc}
2 & 2\\
4 & 0\\
2 & 2\\
2 & 1
\end{array}
\right).
\end{equation}

A group $F$ equipped with an~epimorphism
$F\rightarrow \mathbb{Z}/n\mathbb{Z}$, where $n\in\mathbb{Z}$, is called an~\emph{$n$-indexed} group. This epimorphism $F\rightarrow \mathbb{Z}/n\mathbb{Z}$ is called \emph{degree} and denoted $\deg$. Thus, to any element $f$ of~an~indexed group $F$, an~element $\deg f\in\mathbb{Z}/n\mathbb{Z}$ is assigned, the~group $F$ contains elements of~all degrees and \mbox{$\deg(fg)=\deg f+\deg g$} for any $f$, $g\in F$.

Suppose that $\varphi:F\rightarrow G$ is a~homomorphism from an~$n$-indexed group $F$ to a~group $G$ and $H$ is a~subgroup of~$G$.
The~subgroup
\begin{equation*} 
\begin{array}{l} {
	H_\varphi=\bigcap_{f\in F}H^{\varphi(f)}\cap C(\varphi(\ker\deg))
}\end{array} 
\end{equation*}
is called the~\emph{$\varphi$-core} of~$H$ (see \cite{KM17}). In other words, the~$\varphi$-core $H_\varphi$ of~$H$ consists of~elements $h$ such that $h^{\varphi(f)}\in H$ for all $f$, and $h^{\varphi(f)}=h$ if $\deg f=0$.

\begin{proclaim}{Divisibility theorem \cite{BKV19}}\label{T1}

 Suppose that an~integer $n$ is a~multiple of~the~order of~a~subgroup $H$ of~a~group $G$ and that a~set $\Phi$ of~homomorphisms from an~$n$-indexed group $F$ to $G$ satisfies the~following conditions.
	
	\noindent {\bf I.}
	$\Phi$ is invariant with respect to conjugation by elements of~$H$:
	\begin{center}
		
		if $h\in H$ and $\varphi\in\Phi$, then the~homomorphism
		$\psi \text{: }f\mapsto\varphi(f)^h$ lies in $\Phi$.
		
	\end{center}
	
	\noindent {\bf II.}
	For any $\varphi\in\Phi$ and any element $h$ of~the~$\varphi$-core $H_\varphi$ of~$H$, the~homomorphism $\psi$ defined by
	\begin{equation*} 
		\psi(f)=
		\begin{cases}
			\varphi(f) &\text{for all elements $f\in F$ of~degree zero;}\\
			\varphi(f)h &\text{for some element $f\in F$ of~degree one}\text{ (and, hence, for all degree-one elements)}
		\end{cases}
	\end{equation*}
	belongs to $\Phi$. 
 
 \noindent Then $|\Phi|$ is divisible by $|H|$.
\end{proclaim}
Note that the~mapping $\psi$ from condition {\bf I} is a~homomorphism for any $h\in G$, and the~formula for~$\psi$ from condition {\bf II} defines a~homomorphism for any $h \in H_{\varphi}$ (as explained in \cite{BKV19}).

Let us call tuples of~group elements satisfying a~formula \emph{solutions} (to this formula) henceforth.

A bound variable $t$ is called \emph{isolating} if it occurs in atomic subformulae only in subwords of~the~form $t^{-1}gt$, where $g\in G$. The~coefficient $g$ is called \emph{isolated} in this case.

\begin{remark}\thlabel{r1}
Let $y$ be a~bound variable of~a~formula $\gamma$ over a~group $G$ and $g \in G$. Replacing $y$ with $g^{-1}yg$ or $yg$ in all atomic subformulae of~$\gamma$ does not change the~solution set of~$\gamma$.
\end{remark}
Indeed, for any $g\in G$, if $y$ ranges over the~whole group, then  $y^g$ and $yg$ ranges over the~whole group. That is why two formulae differing only in this are equal: e.g.,  $(\forall y\ \gamma(y,\ z, \dots))\equiv(\forall y\ \gamma(y^h,\ z, \dots))$. 

\begin{lemma}\label{l1}
Let, for an~integer $n$, a~free variable $x_i$ occur in $\gamma$ in the~following forms only:
	
\noindent 1)  $x_i^{kn}$, where $k\in\mathbb{Z}$ ($x_i$ has an~exponent multiple of~$n$);  
	
\noindent 2) $x_ix_jx_i^{-1}$, where $j\neq i$ ($x_i$ conjugates other free variables). 	
 
 \noindent Then the~number of~tuples of~group elements satisfying $\gamma$ is divisible by the~greatest common divisor of~the~centralizer of~a~set of~all non-isolated coefficients of~$\gamma$ and $n$. 
\end{lemma}

\begin{proof} Let $\gamma$ have $m$ free variables $x_1,\dots, x_m$ and put $F=F(x_1, x_2, \dots, x_m)$. Without loss of~generality, assume that $i=1$.  We apply the~divisibility theorem taking $\Phi$ to be a~set of~all homomorphisms  $\varphi\!: F \rightarrow G$ such that tuples  $(\varphi(x_1), \varphi(x_2), \dots, \varphi(x_m))$ satisfy $\gamma$.

As the~indexing deg:\;$F \rightarrow \mathbb{Z}/n\mathbb{Z}$ take an~epimorphism defined by the~mapping of~the~generating set $\text{deg}(x_1)=1$, $\text{deg}(x_2)=\text{deg}(x_3) \dots =\text{deg}(x_m)=\text{deg}(x_1^n)=0$, where $n$ is mentioned in the~first condition. 

Let $H$ be a~subgroup of~the~centralizer of~a~set of~all non-isolated coefficients of~$\gamma$ such that the~order of~$H$ divides $n$. 

To verify that the~first condition of~the~divisibility theorem holds, 
consider a~solution \mbox{$X=(\tilde{x}_1,\  \tilde{x}_2,\dots,\  \tilde{x}_m)$} of~the~formula $\gamma$ and the~corresponding homomorphism \mbox{$\varphi\in \Phi$}, $\varphi(x_i)=\tilde{x}_i$, $i=1\dots m$. Let us check whether the~tuple $(\psi(x_1), \psi(x_2), \dots, \psi(x_m))=(\tilde{x}_1^h,\  \tilde{x}_2^h, \dots, \ \tilde{x}_m^h)$ ($h\in H$) satisfies the~formula $\gamma$.

All atomic subformulae have the~following form: 
\[v_1(y_1)w_1(x_1,\ x_2, \dots,\ x_m)v_2(y_2)w_2(x_1,\ x_2, \dots,\ x_m)\dots v_r(y_r)w_r(x_1,\ x_2, \dots,\ x_m)=1,\]
where $y_i$ are bound variables (not necessarily different), $v_i \in A*\langle y_i\rangle$, $w_i \in A*F(x_1,\ x_2, \dots,\ x_m)$, $A$ is a~subgroup of~the~group $G$ generated by all coefficients of~the~formula $\gamma$.

 We substitute the~new tuple in each atomic subformula keeping in mind that $h$ centralizes the~coefficients of~the~subwords $w_i(x_1,\dots,x_m)$. Then the~left-hand side of~each subformula looks like:
\begin{multline}\nonumber
		v_1(y_1)w_1(\tilde{x}_1^h,\ \tilde{x}_2^h, \dots,\ \tilde{x}_m^h)v_2(y_2)w_2(\tilde{x}_1^h,\  \tilde{x}_2^h, \dots,\ \tilde{x}_m^h)\dots v_r(y_r)w_r(\tilde{x}_1^h, \ \tilde{x}_2^h, \dots,\  \tilde{x}_m^h)=\\
		= v_1(y_1)w_1(\tilde{x}_1,\ \tilde{x}_2,\dots,\  \tilde{x}_m)^hv_2(y_2)w_2(\tilde{x}_1, \ \tilde{x}_2,\dots,\  \tilde{x}_m)^h\dots v_r(y_r)w_r(\tilde{x}_1, \ \tilde{x}_2,\dots, \ \tilde{x}_m)^h=\\
		= v_1(y_1)h^{-1}w_1(\tilde{x}_1,\ \tilde{x}_2,\dots,\  \tilde{x}_m)v_2(y_2)^{h^{-1}}w_2(\tilde{x}_1, \ \tilde{x}_2,\dots,\  \tilde{x}_m)\dots v_r(y_r)^{h^{-1}}w_r(\tilde{x}_1, \ \tilde{x}_2,\dots, \ \tilde{x}_m)h.
	\end{multline}

\noindent Conjugating the~whole atomic subformula by $h^{-1}$, we obtain
\[v_1(y_1)^{h^{-1}}w_1(x_1,\ x_2, \dots,\ x_m)v_2(y_2)^{h^{-1}}w_2(x_1,\ x_2, \dots,\ x_m)\dots v_r(y_r)^{h^{-1}}w_r(x_1,\ x_2, \dots,\ x_m)=1.\]

If the~bound variable $y_i$, $i\in \{1,\dots,r\}$ is non-isolating, then $h$ centralizes the~coefficients of~$v_i$ and hence $v_i(y_i)^{h^{-1}}=v_i(y_i^{h^{-1}})$. Every isolating variable $y_i$ and isolated coefficient $g$ occur in subwords $y_i^{-1}gy_i$ only, but $(y_i^{-1}gy_i)^{h^{-1}}=(y_ih^{-1})^{-1}gy_ih^{-1}$, that is why $v_i(y_i)^{h^{-1}}=v_i(y_ih^{-1})$ in that case. According to the~\thnameref{r1}, such transformations of~the~formula preserve its solution set. Thus, the~first condition of~the~divisibility theorem holds.

To verify the~second condition of~the~divisibility theorem, let us check whether the~tuple \\$(\psi(x_1), \psi(x_2), \dots, \psi(x_m))=(\tilde{x}_1h,\  \tilde{x}_2, \dots, \ \tilde{x}_m)$ is a~solution of~$\gamma$, where $h\in H_{\varphi}$.

\noindent According to the~condition, $x_1$ occurs in $\gamma$ in the~following forms:

\noindent 1) $t^{nk}$, where $k\in\mathbb{Z}$,  
 
\noindent 2) $x_1x_jx_1^{-1}$, where $j\neq 1$. 

Then it suffices to show that $(\tilde{x_1}h)^n=\tilde{x_1}^n$ and $(\tilde{x_1}h)\tilde{x_j}(\tilde{x_1}h)^{-1}=\tilde{x_1}\tilde{x_j}\tilde{x_1}^{-1}$. Since $\tilde{x_j}\in\varphi(\ker\deg)$, where $j\neq 1$, and $h\in H_\varphi\subset C(\varphi(\ker\deg))$ the~second equality is obvious. 

To prove that $(\tilde{x_1}h)^n=\tilde{x_1}^n$, let us use the~following lemma. 
\begin{proclaim}{Brauer lemma \cite{Bra69}}\thlabel{Br} 
	If $U$ is a~finite normal subgroup of~a~group $V$, then, for all $v\in V$, and $u\in U$, the~elements $v^{|U|}$ and $(vu)^{|U|}$ are conjugate by an~element of~$U$.
\end{proclaim}

\noindent Applying the~Brauer lemma to $U=H_\varphi \subset V = H_\varphi \cdot \langle \varphi(x_1)\rangle$, and $v=\varphi(x_1)$, we obtain \mbox{$(\varphi(x_1)h)^{|H_\varphi|}=(\varphi(x_1))^{|H_\varphi|u}$} for some $u\in H_\varphi$. Since $|H_\varphi|$ divides $n$, we have 
\begin{equation}\nonumber
(\tilde{x}_1h)^n=(\varphi(x_1)h)^n=(\varphi(x_1))^{nu}=\varphi(x_1^n)^u=\varphi(x_1^n)=\varphi(x_1)^n=\tilde{x}_1^n.
\end{equation}
Note that deg$(x_1^n)\equiv n \equiv 0$ (mod $n$), and $u\in H_\varphi \subset C(\varphi(\ker\deg))$. This completes the~proof~of~Lemma~1.
\end{proof}
 
\pagebreak

\begin{lemma}\label{l2}
For any formula $\gamma$ with $m$ free variables $x_1,\dots,x_m$ over an~arbitrary group $G$, there exists a~finite formula $\gamma'$ having the~same number of~solutions such that a~free variable $t$ occurs in $\gamma'$ in~the~following forms only:
	
\noindent 1) $t^{nk}$, where $k\in\mathbb{Z}$, $n=\frac{\Delta_m}{\Delta_{m-1}}$,  $\Delta_i$ is the~greatest common divisor of~all minors of~order $i$ of~the~matrix of~the~initial formula $\gamma$, and $\Delta_i=0$, if $i$ is greater than the~number of~atomic subformulae; $\Delta_0=1$; $\frac{0}{0}=0$.  
	
\noindent 2) $t^{-1}zt$, where $z$ is a~free variable of~$\gamma'$ different from $t$. 
 
\end{lemma}

\begin{proof}

\begin{proposition}
	
	Under the~conditions of~the~lemma, the~greatest common divisor of~all coefficients of~the~first column of~the~matrix $A(\gamma)$ is divisible by $\frac{\Delta_m}{\Delta_{m-1}}$ after a~suitable invertible change of~free variables. 
\end{proposition}

\begin{proof}It is well known that using invertible integer transformations of~columns and rows, the~matrix $A(\gamma)$ can be diagonalized in such a way that each diagonal element divides the~next element. The~integer $\frac{\Delta_m}{\Delta_{m-1}}$ is the~$m$th invariant factor of~the~matrix of~the~formula. By forgetting about the~row transformations and performing the~column transformations only, we obtain a~matrix of~arbitrary form, but the~elements of~the~$m$th column are still divisible by $\frac{\Delta_m}{\Delta_{m-1}}$. Elementary column transformations are induced by invertible changes of~variables $x_i \rightarrow x_ix_j^k$. The~transformation adds the~$i$th column multiplied by $k$ to the~$j$th column. This completes the~proof.
\end{proof}

Thus, we can assume that 
\begin{center}
\begin{minipage}[c][0.5cm][c]{0.92\textwidth}
		\emph{the~exponent sum of~the~free variable $t$ (for example, $x_1$) is a~multiple of~$n = \frac{\Delta_m}{\Delta_{m-1}}$ in each cycle of~the~graph of~the~formula (and in each atomic subformula).}
\end{minipage}\hfill$(*)$
\end{center}
Let us fix some connected component of the~graph  $\Gamma(\gamma)$ and consider an~arbitrary path  $p_{ij}$ from a~vertex  $y_i$ to a~vertex $y_j$ in this connected component. Property $(*)$ implies that the~sum of~the~first coordinates of~the~labels of~the~edges of~$p_{ij}$ modulo $n$ does not depend on the~choice of~the~path. Denote by $s_{ij}$ the~residue of this sum modulo $n$.

Next, we proceed as follows.

\noindent\textbf{Step 1.} Property $(*)$ implies the~feasibility of~bringing together the~exponent of~the~free variable $t$ to multiples of~$n$ by conjugating coefficients, free variables, and bound variables of~the~formula with $t$ to some exponent $0\leq k < n$. Namely, we rewrite all atomic subformulae as equations $w' = 1$, where $w'$ is a~word depending on subwords of~the~form $t^n$, $ t^{-k}zt^{k} $, and $ t^{-k}(y^{-1}gy)t^{k} $, where $z$ is either a~free variable ($z \neq t$), or a~non-isolated coefficient, or a~non-isolating bound variable, $y$ denotes an~isolating bound variable, $g$ stands for an~isolated coefficient, and the~exponent $0\leq k < n$ may vary for different variables and different occurrences of~the~same variable.

For example, consider the~formula $\gamma(t, x)$:
\begin{equation*}
\exists y_1, y_2, y_3, y_4 (y_1txy_3t^2x^2y_1t^3x^3=1) \wedge (y_3tgy_2tgy_1t=1) \lor \neg(h^{y_4}t^5y_3tx=1), 
\end{equation*}
where $h$ and $g$ are some elements of~the~group $G$.

We have the~matrix 
$A(\gamma) = \left(
\begin{array}{cc}
3&3 \\
3&1\\
3&0\\
0&2\\
6&1
\end{array}
\right)$
and the~graph $\Gamma(\gamma)$:

\begin{figure}[h!]
    \centering
    \includegraphics[width=0.6\textwidth]{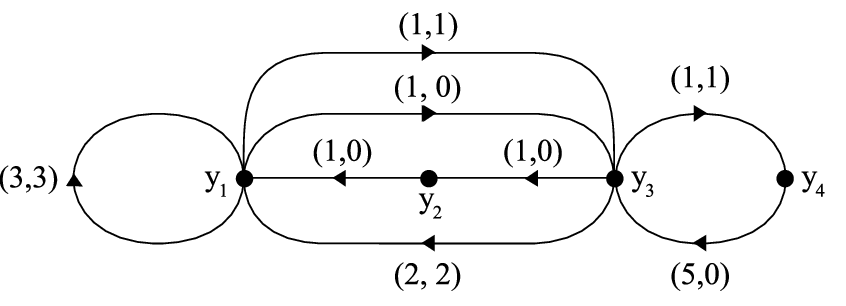}
    \caption{The~graph of~the~formula $\gamma(t)$}
\end{figure}

We obtain $n = \frac{\Delta_2}{\Delta_1} = \frac{3}{1} = 3$. As a~result of~Step~1, the~atomic subformulae of~$\gamma(t, x)$ take the~following form:

\[
\begin{aligned}
& y_1t^3x^{t^2}y_3^{t^2}x^2y_1t^3x^3 = 1,\\
& y_3t^3g^{t^2}y_2^{t^2}g^ty_1^t = 1,\\
& h^{y_4}t^6y_3^tx = 1.
\end{aligned}
\]

Notice that property $(*)$ implies a~specific dependence between different occurrences of~bound variables:
\begin{center}
\begin{minipage}[c][0.5cm][c]{0.92\textwidth}
\emph{If $y_i^{t^{k}}$ (or $a^{y_it^{k}}$) and $y_j^{t^{l}}$ (or $ b^{y_jt^{l}}$) occur in the~atomic subformula $w'$, then $k - l \equiv s_{ij}$ (mod $n$), where $y_i$ and $y_j$ are bound variables, $a$ and $b$ are isolated coefficients of~the~formula $\gamma$}. 
\end{minipage}\hfill $(**)$
\end{center}

In particular, this yields that all occurrences of~a~bound variable in $w'$ are \textit{consistent}, namely, if $i = j$, then $k = l$.

\noindent\textbf{Step 2.} We want to ensure consistency of~occurrences of~an~arbitrary bound variable $y_q$ in all atomic subformulae. If $y_q^{t^{k}}$ (or $g^{y_qt^{k}}$) occurs in the~equation $w' = 1$, we replace it with the~equivalent equation $(w')^{t^{-k}} = 1$. We transform each atomic subformula containing $y_q$ in this manner.

Notice that property $(**)$ implies that Step 2 leads to the~consistency of~all occurrences of~selected atomic subformulae for each bound variable. Indeed, as a~result of~Step 2, the~fixed bound variable $y_i$ occurs in these atomic subformulae only in subwords of~the~form $y_i^{t^{s}}$ (or $g^{y_it^{s}}$), where $s \equiv s_{iq}$ (mod $n$).

For the~formula $\gamma(t, x)$ from the~example, choose $y_q = y_1$. In this case, only the~second equation needs to be modified. By conjugating both sides by $t^{-1}$, the~atomic subformula takes the~form $t^3y_3^{t^2}g^ty_2^tgy_1 = 1 $.

\noindent\textbf{Step 3.} Similarly, we want to ensure the~consistency of~all occurrences of~bound variables in atomic subformulae $w'$ not containing $y_q$, replacing them with conjugates using $t^{s_{iq} - k}$, where $i$ is the~index of~the~bound variable and $k$ is the~corresponding exponent of~$t$ in the~given subformula.

After this adjustment, any other bound variable $y_j$ occurs in $(w')^{t^{s_{iq}-k}}$ only as subwords of~the~form $y_j^{t^{s_{jq}}}$ or $g^{y_jt^{s_{jq}}}$. Indeed, $l + s_{iq} - k \equiv s_{ji} + s_{iq} \equiv s_{jq}$ (mod $n$), where $l$ is the~degree of~the~free variable $t$ in the~subwords $y_j^{t^l}$ (or $g^{y_jt^l}$) of~$w'$.

Thus, we have ensured the~consistency of~each bound variable in all occurrences.

For example, after Step 3, the~third atomic subformula of~$\gamma(t, x)$ is modified, and all occurrences of~bound variables are consistent:

\[
\begin{aligned}
& y_1t^3x^{t^2}y_3^{t^2}x^2y_1t^3x^3 = 1,\\
& t^3y_3^{t^2}g^ty_2^tgy_1 = 1,\\
& h^{y_4t}t^6y_3^{t^2}x^t = 1.
\end{aligned}
\]

\noindent\textbf{Step 4.} Now, each bound variable $y$ occurs in the~formula only in the~subwords of~the~form $y^{t^{k}}$ or $yt^{k}$, where $k$ is the~same for all occurrences of~the~variable $y$. Therefore, we can discard all conjugations or multiplications by $t^{k}$ everywhere because of~the~\thnameref{r1}.

After Step 4 the formula $\gamma(t, x)$ looks as follows: \[
\exists y_1, y_2, y_3, y_4 (y_1t^3x^{t^2}y_3x^2y_1t^3x^3 = 1)\wedge
(t^3y_3g^ty_2gy_1 = 1)\lor \neg (h^{y_4}t^6y_3x^t = 1).
\]

Let us do Steps 1-4 for each connected component of the graph
$\Gamma(\gamma)$. Since the component can be
chosen arbitrarily, the conclusion of Step 4 holds for any variable.

\noindent\textbf{Step 5.} For each free variable and non-isolated coefficient, we introduce new symbols $z^{(i)}$, where \mbox{$0 \leqslant i < n$}, and replace the~obtained formula as follows:

\[\gamma'=
\hat{\gamma} \wedge \underbrace{ \bigwedge_{z}\left(\left(z^{(0)}=z\right) \wedge \left(z^{(1)}=t^{-1}z^{(0)}t\right) \wedge \dots \wedge \left(z^{(n-1)}=t^{-1}z^{(n-2)}t\right)\right)}_{\textstyle \gamma_1},
\]

\noindent where $z$ ranges over all free variables other than $t$ and non-isolated coefficients of~the~initial formula $\gamma$, and $\hat{\gamma}$ is obtained from the~transformed $\gamma$ by replacing all subwords $z^{t^i}$ with $z^{(i)}$. We assume symbols $z^{(i)}$ to be free variables, $0 \leqslant i < n $.

For the~formula $\gamma(t, x)$ from the~example, we obtain $\gamma'(t, x, x^{(0)}, x^{(1)}, x^{(2)})$:
\begin{multline*}\exists y_1, y_2, y_3, y_4
\bigg(\left(y_1t^3x^{(2)}y_3(x^{(0)})^2y_1t^3(x^{(0)})^3 = 1\right)\wedge
\left(t^3y_3g^{(1)}y_2g^{(0)}y_1 = 1\right)\lor \neg
\left(h^{y_4}t^6y_3x^{(1)} = 1\right)\bigg)\wedge\\\wedge \bigg( \left(g^{(0)}=g\right) \wedge \left(g^{(1)}=t^{-1}g^{(0)}t\right) \wedge\ \left(x^{(0)}=x\right) \wedge\left(x^{(1)}=t^{-1}x^{(0)}t\right) \wedge\left(x^{(2)}=t^{-1}x^{(1)}t\right) \bigg).
\end{multline*}
Obviously, the~numbers of~tuples satisfying the~new formula $\gamma'$ and the~initial formula $\gamma$ are equal. 

The~proof~is provided for the~case $n > 0$. In the~case $n = 0$, we remove the~upper limits on the~indices, and all congruences modulo become exact equations. 

According to the~transformations made, the~free variable $t$ occurs in $\hat{\gamma}$ in the~form $t^{nk}$ only, where $k\in\mathbb{Z}$. Indeed, we have brought together the~exponent of~$t$ in Step 1. In Steps 2-3, we have conjugated subformulae by corresponding exponents of~$t$. In the subformula $\gamma_1$, $t$ occurs in the form $t^{-1}zt$ only, where $z$ is a~free variable different from $t$. Thus, the~formula $\gamma'$ satisfies the~conditions of~the~lemma. 
\end{proof}

\begin{proclaim}{Main theorem}
	In any group, the~number of~tuples of~group elements satisfying a~formula $\gamma$ with $m$ free variables is divisible by the~greatest common divisor of~the~centralizer of~the~set of~all non-isolated coefficients of~$\gamma$ and $n=\frac{\Delta_m}{\Delta_{m-1}}$, where $\Delta_i$ is the~greatest common divisor of~all minors of~order $i$ of~the~matrix of~the~formula. The~following conventions are assumed: $\Delta_i=0$ if $i$ is larger than the~number of~rows of~the~matrix; $\Delta_0=1$; $\frac{0}{0}=0$.
	
\end{proclaim}

\begin{proof} Lemma~\ref{l1} and Lemma~\ref{l2} easily imply the~main theorem. Indeed, according to Lemma~\ref{l2}, for any formula $\gamma$ there exists a~formula $\gamma'$ satisfying conditions 1) and 2) of~Lemma~\ref{l1}  with the~same number of~solutions, and  $n=\frac{\Delta_m}{\Delta_{m-1}}$, where $\Delta_i$ is the~greatest common divisor of~all minors of~order $i$ of~the~matrix of~the~initial formula. According to Lemma~\ref{l1}, the~number of~solutions to $\gamma'$ in $G$ (hence, the~number of~solutions to $\gamma$) is divisible by the~greatest common divisor of~the~centralizer of~the~set of~all non-isolated coefficients of~$\gamma$ and $n$.
\end{proof}

\section{Examples}
Our research began with the~Frobenius theorem, therefore, the~most basic examples are first-order formulae with one free variable and nonzero exponent sum in at least one atomic subformula. For instance, the~number of~squares with cube root in a~group $G$ is divisible by GCD$(G,\;2)$, and the~number of~cubes with square root in the~group $G$ is divisible by GCD$(G,\;3)$. 

We expand the~set of~examples by adding some conditions to the~formulae. For any given integers $n$ and $k$, GCD$(G,\;n)$ divides the~number of~elements satisfying formulae 1-3 with one free variable, and GCD$(G,\;2n)$ divides the~number of~elements satisfying formulae 4-6 with two free variables:
\begin{enumerate}
\item $\exists z \;(x^n=z^k)\wedge\neg(x^n=1)$;
\item $\exists z,\;t \;(x^n=z^k)\wedge\neg([z,\;t]=1)$;
\item $\exists z,\;t \;(x^n=[z,\;t])\wedge\neg(x^n=1)$;
\item $\exists z \;(x^ny^n=z^k)\wedge\neg(x=y)$;
\item $\exists z \;((xy)^n=z^k)\wedge\neg(x^n=y^n)$;
\item $\exists z,\;t\;((xy)^n=[t,\;z]^k)\wedge(y^nx^n=[z^k,\;t^k])\wedge\neg(x=y)$.
\end{enumerate}
\noindent the~matrixes are
$
A(\gamma_1) = \left(
\begin{array}{cc}
n\\
n
\end{array}
\right),\;  
A(\gamma_2) = \left(
\begin{array}{cc}
n\\
0\\
0\\
0
\end{array}
\right),\; 
A(\gamma_3) = \left(
\begin{array}{cc}
n\\
0\\
0\\
n
\end{array}\right),\;
A(\gamma_4) = \left(
\begin{array}{cc}
n&n\\
1&-1
\end{array}
\right),\; 
$

\noindent $
A(\gamma_5) = \left(
\begin{array}{cc}
n&n\\
n&-n
\end{array}
\right),\; 
A(\gamma_6) = \left(
\begin{array}{cc}
n&n\\
n&n\\
0&0\\
\vdots&\vdots\\
0&0\\
1&-1
\end{array}
\right). 
$

\noindent The~graphs of~the~formulae look as follows:

\begin{figure}[ht]
    \centering
    \subfigure(1) {\includegraphics[width=0.29\textwidth]{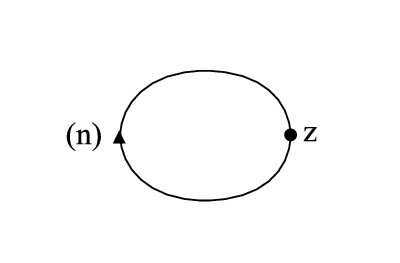}}
    \subfigure(2){\includegraphics[width=0.29\textwidth]{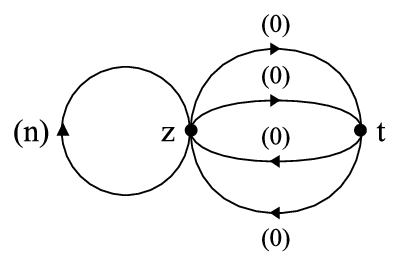}} 
    \subfigure(3) {\includegraphics[width=0.29\textwidth]{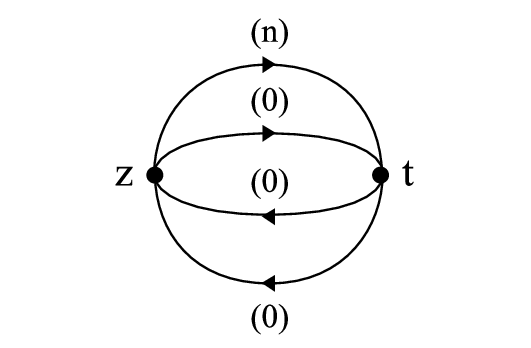}}
    \caption{Graphs for examples 1-3.}
\end{figure}

\begin{figure}[ht]
    \centering
    \subfigure(4) {\includegraphics[width=0.29\textwidth]{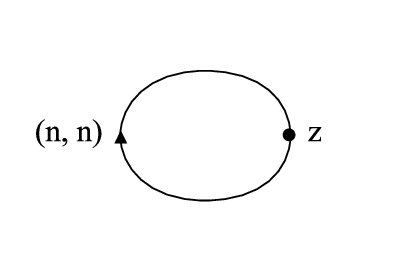}}
    \subfigure(5){\includegraphics[width=0.29\textwidth]{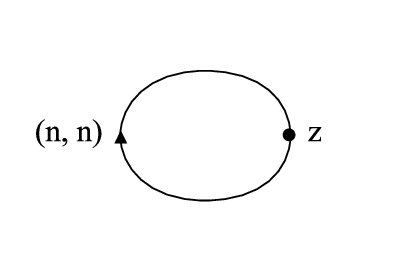}}
    \subfigure(6) {\includegraphics[width=0.29\textwidth]{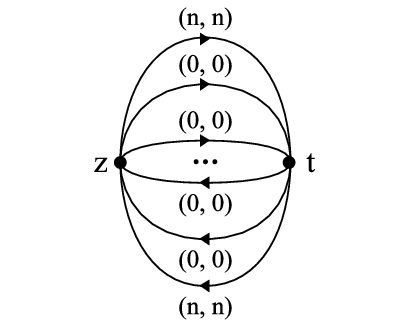}}
    \caption{Graphs for examples 4-6.}
\end{figure}


\newpage

\begin{thebibliography}{22}
	
	
	
	
	\bibitem{ACNT13}
	T. Asai, N. Chigira, T. Niwasaki, Yu. Takegahara, 
	On a~theorem of~P. Hall,
	Journal of~Group Theory, 16:1 (2013), 69-80.
 
	\bibitem{AsTa01}
	T. Asai, Yu. Takegahara,
	$|\text{Hom}(A,G)|$, 
	IV, J. Algebra, 246 (2001), 543-563.
	
	\bibitem{Bra69}
	R. Brauer, 
	On a~Theorem of~Frobenius,
	The American Mathematical Monthly, 76:1 (1969), 12-15.
	\bibitem{BrTh88}
	K. Brown, J. Th\'evenaz,
	A generalization of~Sylow's third theorem,
	J. Algebra, 115 (1988), 414-430.
        \bibitem{BKV19}
	E. K. Brusyanskaya, A. A.Klyachko, A. V. Vasil'ev,
	What do Frobenius's, Solomon's, and Iwasaki's theorems on divisibility in groups have in common?,
	Pacific Journal of~Mathematics, 302:2 (2019), 437-452. See also arXiv:1806.08870.

 \bibitem{BK22}E. K. Brusyanskaya, A. A. Klyachko, On the~number of~epi-, mono- and homomorphisms of~groups, Izv. Math., 86:2 (2022), 243–251. See also arXiv:2012.03123.
 
	\bibitem{Frob95}
	F. G. Frobenius,
	Verallgemeinerung des Sylow'schen Satzes, Sitzungsberichte der K\"onigl. 
	Preu. Akad. der Wissenschaften (Berlin) (1895), 981-993.
	\bibitem{Frob03}
	F. G. Frobenius, 
	\"Uber einen Fundamentalsatz der Gruppentheorie, Sitzungsberichte der K\"onigl. 
	Preu. Akad. der Wissenschaften (Berlin) (1903), 987-991.
	\bibitem{GRV12}
	C. Gordon, F. Rodriguez-Villegas,
	On the~divisibility of~$\#\text{Hom}(\Gamma, G)$ by $|G|$, 
	J. Algebra, 350:1 (2012), \text{300-307}. See also arXiv:1105.6066.
	\bibitem{Hall36}
	Ph. Hall, 
	On a~theorem of~Frobenius, 
	Proc. London Math. Soc. 40 (1936), 468-531.
	\bibitem{Isaa70}
	I. M. Isaacs,
	Systems of~equations and generalized characters in groups,
	Canad. J. Math., 22 (1970), \hbox{1040-1046}.
	
	\bibitem{KM14}
	A. A. Klyachko, A. A. Mkrtchyan,
	How many tuples of~group elements have a~given property? With an~appendix by Dmitrii V. Trushin,
	Intern. J. of~Algebra and Comp. 24:4 (2014), 413-428. See also arXiv:1205.2824.
	\bibitem{KM17}
	A. A. Klyachko, A. A. Mkrtchyan,
	Strange divisibility in groups and rings, 
	Arch. Math. 108:5 (2017), 441-451. See also arXiv:1506.08967.
	\bibitem{Kula38}
	A. Kulakoff, Einige Bemerkungen zur Arbeit: “On a~theorem of~Frobenius” von P. Hall, Math. Sb. 3:2 (1938), 403–405.
	
	\bibitem{Sehg62}
	S. K. Sehgal, On P.
	Hall's generalisation of~a~theorem of~Frobenius,
	Proc. Glasgow Math. Assoc.,  5 (1962), 97-100.
	\bibitem{Solo69}
	L. Solomon, 
	The solutions of~equations in groups, 
	Arch. Math., 20:3 (1969), 241--247.
	\bibitem{Yosh93}
	T. Yoshida,
	$|\text{Hom}(A, G)|$, 
	Journal of~Algebra, 156:1 (1993), 125-156.
	
\end{thebibliography}
\end{document}